\documentclass[12pt]{amsart} 
\usepackage{amssymb,colordvi,multirow,graphicx,amsmath} 

\newtheorem{thm}{Theorem}

\newtheorem{lem}[thm]{Lemma}

\newtheorem{example}[thm]{Example}
\newtheorem{remark}[thm]{Remark}

\newcounter{FNC}[page]
\def\fauxfootnote#1{{\addtocounter{FNC}{2}$^\fnsymbol{FNC}$%
     \let\thefootnote\relax\footnotetext{$^\fnsymbol{FNC}$#1}}}

 \copyrightinfo{}{}

\newcommand{\eps}{{\varepsilon}}
\newcommand{\R}{\mathbb{R}}
\newcommand{\Log}{\mbox{Log}}
\newlength{\khov}
\renewcommand{\qed}{$\blacksquare$}

\title[Sharpness of fewnomial bounds]{On the Sharpness of fewnomial bound and 
  the number of components of a fewnomial hypersurface} 

\author{Fr\'ed\'eric Bihan}
\address{Laboratoire de Math\'ematiques\\
         Universit\'e de Savoie\\
         73376 Le Bourget-du-Lac Cedex\\
         France}

\email{Frederic.Bihan@univ-savoie.fr}
\urladdr{http://www.lama.univ-savoie.fr/\~{}bihan/}

\author{J.~Maurice Rojas}
\address{Department of Mathematics\\
         Texas A\&M University\\
         College Station\\
         Texas \ 77843\\
         USA}
\email{rojas@math.tamu.edu}
\urladdr{http://www.math.tamu.edu/\~{}rojas/}
\thanks{Rojas supported by NSF CAREER grant DMS-0349309}  

\author{Frank Sottile}
\address{Department of Mathematics\\
         Texas A\&M University\\
         College Station\\
         Texas \ 77843\\
         USA}
\email{sottile@math.tamu.edu}
\urladdr{http://www.math.tamu.edu/\~{}sottile/}
\thanks{Sottile supported by NSF CAREER grant DMS-0538734}  

\keywords{Fewnomials, connected component}
\subjclass[2000]{14P99}

\begin{document}

\begin{abstract}
 We prove the existence of systems of $n$ polynomial equations in $n$ 
 variables with a total of $n+k+1$ distinct 
 monomial terms possessing $\left\lfloor \frac{n+k}{k}\right\rfloor^k$ nondegenerate
 positive solutions.
 This shows that the recent upper 
 bound of  $\frac{e^2+3}{4}\, 2^{\binom{k}{2}} n^k$ for the number of 
 nondegenerate positive solutions has the correct order for 
 fixed $k$ and large $n$.
 We also adapt a method of Perrucci to show that there are fewer than
 $\frac{e^2+3}{4}2^{\binom{k+1}{2}}2^nn^{k+1}$ 
connected components in a smooth hypersurface
 in the positive orthant of $\R^N$ defined by a polynomial with $n{+}k{+}1$
 monomials, where $n$ is the dimension of the affine span of the 
exponent vectors. 
 Our results hold for polynomials with real exponents.
\end{abstract}
\maketitle

\section{Introduction} 
Khovanskii~\cite{Kh80} showed there are at most 
$2^{\binom{n+k}{2}}(n+1)^{n+k}$ nondegenerate positive solutions to a 
system of $n$ polynomial equations in $n$ variables  
which are linear combinations of (the same) $n{+}k{+}1$ monomials.
This  fewnomial bound is also valid in the more general 
setting of linear combinations of monomials with real-number
exponents. 
The underlying bounds are identical whether one uses 
integral or real exponents~\cite{Na01}, and the arguments 
of~\cite{BS} require that we allow real exponents. 

While Khovanskii's fewnomial bound was not believed to be sharp, only 
recently have
smaller bounds been found. 
The first breakthrough was due to Li, Rojas, and Wang~\cite{LRW03} who
showed that a system of two trinomials in two variables has at most 5 positive
solutions --- 
%
%
which is
smaller than Khovanskii's bound of $5184$.
Bihan~\cite{Bihan} showed that a system of $n$ polynomials in $n$ variables
with $n{+}2$ monomials has at most $n{+}1$ nondegenerate positive solutions 
and proved the existence of such a system with $n{+}1$ positive solutions.
Bihan and Sottile~\cite{BS} generalized this to all $k$, giving
the upper bound of $\frac{e^2+3}{4}\, 2^{\binom{k}{2}} n^k$ for the number
of nondegenerate positive solutions, which is significantly smaller than
Khovanskii's bound.

\subsection{A Lower Bound for Fewnomial Systems} 
We show that the  Bihan-Sottile upper bound~\cite{BS} is near-optimal 
for fixed $k$ and large $n$.

\begin{thm}\label{Th:one}
  For any positive integers $n,k$ with $n>k$, there exists a system of $n$ polynomials in 
  $n$ variables involving $n{+}k{+}1$ distinct monomials 
  and having $\left\lfloor \frac{n+k}{k}\right\rfloor^k$ nondegenerate positive solutions. 
\end{thm} 

\noindent
We believe that there is room for improvement in the dependence on $k$, both in the upper
bound of~\cite{BS} and in the lower bound of Theorem~\ref{Th:one}.

\begin{proof}
 We will construct such a system when $n=km$, a multiple of $k$, from which 
we may deduce the
 general case as follows.
 Suppose that $n=km+j$ with $1\leq j<k$ and we have a system of $mk$ equations in $mk$ variables 
 involving $mk{+}k{+}1$ monomials with $(m{+}1)^k$ nondegenerate positive solutions.
 We add  $j$ new variables $x_1,\ldots,x_j$ and $j$ new equations $x_1=1, \ldots,x_j=1$.
 Since the polynomials in the original system may be assumed to have constant 
terms, this gives a system with $n$ polynomials in $n$ variables having 
$n+k+1$ monomials and $(m{+}1)^k=\left\lfloor \frac{n+k}{k}\right\rfloor^k$ 
nondegenerate positive solutions. So let us fix positive integers $k,m$ and 
set $n=km$.\smallskip 

 Bihan~\cite{Bihan} showed there exists a system of $m$ polynomials in $m$
 variables 
\[
   f_1(y_1,\ldots,y_m)\ =\ \cdots \ =\ 
   f_m(y_1,\ldots,y_m)\ =\ 0
\]
 having $m{+}1$ solutions, and where each polynomial has the same 
 $m{+}2$ monomials, one of which we may take to be a constant.

 For each $j=1,\ldots,k$, let $y_{j,1},\ldots,y_{j,m}$ be $m$ new variables
 and consider the system
 \[
   f_1(y_{j,1},\ldots,y_{j,m})\ =\ \cdots \ =\ 
   f_m(y_{j,1},\ldots,y_{j,m})\ =\ 0\,,
 \]
 which has $m{+}1$ positive solutions in
 $(y_{j,1},\ldots,y_{j,m})$. 
 As the sets of variables are disjoint, the combined
 system consisting of all $km$ polynomials in all $km$ variables has $(m{+}1)^k$ 
 positive solutions. 
 Each subsystem has $m{+}2$ monomials, one of which is a constant.
 Thus the combined system
 has $1+k(m{+}1)=km{+}k{+}1=n{+}k{+}1$ monomials.
\end{proof} 

\begin{remark} {\rm 
 Our proof of Theorem~\ref{Th:one} is based on Bihan's 
 nonconstructive proof 
 of the existence of a system of $n$ polynomials in $n$ variables having 
 $n{+}2$ monomials
 and $n{+}1$ nondegenerate positive solutions.
 While finding such systems explicitly is challenging in general, let us 
 do so for $n\in\{2,3\}$.

 The system of $n=2$ equations with $2$ variables 
\[
    x^2y - (1+4x^2)\ =\  xy - (4+x^2)\ =\ 0\,,
\]
 has $4=2+1+1$ monomials and exactly $3$ complex solutions, each of which is nondegenerate and 
 lies in the positive quadrant. 
 We give numerical approximations, computed with the computer algebra system
 SINGULAR~\cite{GPS05,PW},
\[
   (2.618034,\, 4.1459), \quad (1,\, 5), \quad 
     (0.381966,\, 10.854102)\,.
\]

 The system of $n=3$ equations with $3$ variables 
\[
   yz - (1/8+2x^2)\ =\  xy - (1/220 + x^2)\ =\   z - (1+x^2)\ =\  0\,,
\]
 has $5=3+1+1$ monomials and exactly $4$ complex solutions, each of which is 
nondegenerate and  lies in the positive octant,
 \begin{eqnarray*}
  &   (0.076645,\, 0.1359,\, 1.00587), \quad 
   (0.084513,\, 0.13829,\, 1.00714)\,,&\\
&(0.54046,\, 0.54887,\, 1.2921),\quad 
   (1.29838,\, 1.30188,\, 2.6858)\,.& 
\end{eqnarray*}}
\end{remark}

\subsection{An Upper Bound for Fewnomial Hypersurfaces} 
Khovanskii also considered smooth hypersurfaces in the positive orthant 
$\R^n_>$ defined by polynomials with $n{+}k{+}1$ monomials.
He showed~\cite[Sec.\ 3.14, Cor.~4]{Kh91} that the total Betti number
 of such a fewnomial hypersurface  is at most
\[
(2n^2-n+1)^{n+k}(2n)^{n-1}2^{\binom{n+k}{2}}\,.
\]
Li, Rojas, and Wang \cite{LRW03} bounded the number of
connected components of such a hypersurface by 
$n(n+1)^{n+k+1}2^{n-1}2^{\binom{n+k+1}{2}}$.
(This bound holds also for singular hypersurfaces.)
Perrucci~\cite{Pe05} lowered this bound to $(n+1)^{n+k}2^{1+
\binom{n+k}{2}}$. 
His method was to bound the number of compact components and then use an 
argument based on the faces of the $n$-dimensional cube to bound the number 
of all components.  We improve Perrucci's method, using the $n$-simplex and 
the bounds of Bihan and Sottile~\cite{BS} to obtain a new, lower bound.

\begin{thm}\label{Th:two}
  A smooth hypersurface in $\R^N_>$ defined by a polynomial with $n{+}k{+}1$ monomials
  whose exponent vectors have $n$-dimensional affine span has fewer than
\[
   \frac{e^2+3}{4}\cdot 2^{\binom{k+1}{2}}2^nn^{k+1}
\]
 connected components. 
\end{thm}

Both our method and that of Perucci estimate the number of compact components in the
intersection of the hypersurface with linear spaces which are chosen to preserve the
sparse structure so that the fewnomial bound still applies.
 For this same reason, this same method was used by Benedetti, Loeser, and
Risler~\cite{blr}  to bound the number of components of a real algebraic set in terms of
the degrees of its defining polynomials.

\begin{remark}{\rm 
  If a hypersurface $f=0$ is singular, we may bound its number of connected components by the
  number of connected components of the two hypersurfaces $f=\pm\eps$, for 
 $\eps>0$ small. 
 (This is Lemma~13 in~\cite{LRW03}.)
 Since we can assume that $f$ has a constant term, this gives the bound of 
 $2\cdot \frac{e^2+3}{4}\cdot 2^{\binom{k+1}{2}}2^nn^{k+1}$ 
for the number of connected 
 components of an arbitrary (possibly singular) hypersurface in $\R^N_>$ defined by a
 polynomial with $n{+}k{+}1$ monomials whose exponent vectors have an $n$-dimensional
  affine span.  
}\end{remark}

Observe that if the exponent vectors of the monomials in a polynomial $f$ in $N$ variables
have $n$-dimensional affine span, then there is a monomial change of variables on $\R^N_>$
so that $f$ becomes a polynomial in the first $n$ variables, and thus our hypersurface
becomes a cylinder 
\[
   \R^{N-n}_<\ \times\ \{x\in\R^n_>\mid f(x)=0\}\,,
\]
which shows that it suffices to prove Theorem~\ref{Th:two} when $n=N$.
That is, for smooth hypersurfacs in $\R^n$ defined by polynomials with $n{+}k{+}1$
monomial terms whose 
exponent vectors affinely span $\R^n$.

Let $\kappa(n,k)$ be the maximum number of compact connected components of
such a smooth hypersurface and let $\tau(n,k)$ be the maximal number of
connected components of such a hypersurface. 
We deduce Theorem~\ref{Th:two} from the following combinatorial lemma.

\begin{lem}\label{L:three}
 \ \ ${\displaystyle 
   \tau(n,k)\ \leq\ \sum_{i=0}^{n-1} \binom{n+1}{i} \kappa(n-i,k+1)\,}$.
\end{lem}

{\it Proof of Theorem~$\ref{Th:two}$}.
Bihan and Sottile~\cite{BS} proved that 
\[
\kappa(n,k) \leq \frac{e^2+3}{8}2^{\binom{k}{2}}n^k\,.
\]
Substituting this into Lemma~\ref{L:three} bounds $\tau(n,k)$ by
\begin{multline*}
\frac{e^2+3}{8}2^{\binom{k+1}{2}}  \sum_{i=0}^{n-1} \binom{n+1}{i} (n-i)^{k+1}
       \\
   <\  \frac{e^2+3}{8}2^{\binom{k+1}{2}}
             n^{k+1}\sum_{i=0}^{n+1} \binom{n+1}{i}
        \ =\ \frac{e^2+3}{8}2^{\binom{k+1}{2}} n^{k+1} 2^{n+1}\,. \qquad\mbox{\qed}
\end{multline*}


 {\it Proof of Lemma~$\ref{L:three}$.}
 Let $f$ be a polynomial in the variables $x_1,\ldots,x_n$ which has $n{+}k{+}1$
 distinct monomials whose exponent vectors affinely span $\R^n$ and suppose that $f(x)=0$
 defines a smooth hypersurface $X$ in $\R^n_>$. 
 We may apply a monomial change of coordinates to $\R^n_>$ and assume that 
 $1,x_1,x_2,\ldots,x_n$ are among the monomials appearing in $f$.

 Suppose that $\eps:=(\eps_0,\eps_1,\ldots,\eps_n)\in\R^{1+n}_>$ 
 with $\eps_0 \eps_1\cdots\eps_n\neq 1$.
 Define hypersurfaces $H_0,H_1,\ldots,H_n$ of $\R^n_>$ by
\begin{eqnarray*}
   H_0&:=& \{x\in\R^n_>\mid x_1\dotsb x_n\eps_0=1\}\,\quad\mbox{and} \\
   H_i& :=& \{x\mid x_i=\eps_i\}\quad \mbox{for}\quad i=1,\ldots,n\,.
\end{eqnarray*}
 The transformation $\Log\colon\R^n_>\to\R^n$ defined on each coordinate by 
 $x_i\mapsto \log(x_i)$ sends the hypersurfaces $H_0,H_1,\ldots,H_n$ to
 hyperplanes in general position.
 That is, if $S\subset\{0,1,\ldots,n\}$ and we define $H_S:=\cap_{i\in S}H_i$,
 then $\Log(H_S)$ is an affine linear space of dimension $n{-}|S|$.
 Moreover, the complement of the union of hypersurfaces $H_i$ has $2^{n{+}1}-1$ connected
 components, exactly one of which is bounded.

 If we restrict $f$ to some $H_S$, we will obtain a new polynomial $f_S$ in
 $n-|S|$ variables with at most $1{+}(n{-}|S|){+}(k{+}1)$ monomials.
 Indeed, if $i\in S$ with $i\neq 0$, then the equation $x_i=\eps_i$ allows us to eliminate
 both the variable  and the monomial $x_i$ from $f$.
 If however $0\in S$, then we pick an index
 $j\not\in S$ and use $x_1\dotsb x_n\eps_0=1$ to eliminate the variable $x_j$,
 which will not necessarily eliminate a monomial from $f$.
 For almost all $\eps$, the polynomial $f_S$ defines a smooth hypersurface
 $X_S$ of $H_S$.

 We may choose $\eps$ small enough so that every compact connected component
 of $X$ lies in the bounded region of the complement of the hypersurfaces $H_i$,
 and every noncompact connected component of $X$ meets some hypersurface $H_i$.
 Shrinking $\eps$ if necessary, we can ensure that 
 every bounded component of $X_S$ lies in the bounded region of the complement
 of $H_j\cap H_S$ for $j\not\in S$, and every unbounded component meets some 
 $H_j\cap H_S$ for $j\not\in S$.

 Given a connected component $C$ of $X$, the subsets $S\subset\{0,1,\ldots,n\}$
 such that $C$ meets $H_S$ form a simplicial complex.
 If $S$ represents a maximal simplex in this complex, then $C\cap H_S$ is a
 union of compact components of $X_S$, and $|S|<n$ as $H_S$ is not a point.
 Thus the number of connected components of $X$ is bounded by the sum of the
 numbers of compact components of $X_S$ for all $S\subset\{0,1,\ldots,n\}$ with
 $n>|S|$. 
 Since each $f_S$ has at most $1{+}(n{-}|S|){+}(k{+}1)$ monomials,
 this sum is bounded by the sum in the statement of the lemma.
\endproof

\begin{remark}{\rm 
  If $f$ contains a monomial $x^a:=x_1^{a_1}x_2^{a_2}\cdots x_n^{a_n}$ with 
  no $a_i=0$, then we can alter the proof of Lemma~\ref{L:three}
  to obtain a bound of
\[
   \frac{e^2+3}{4}\cdot 2^{\binom{k}{2}}2^nn^{k}
\]
 connected components for the hypersurface defined by $f$.
 
 The basic idea is that if we redefine $H_0$ to be 
\[
   H_0\ :=\ \{x\in\R^n_>\mid x^a\eps_0=1\}\,.
\]
 then the polynomials $f_S$ on $H_S$ have only
  $1+(n{-}|S|)+k$ monomials, and so we estimate the number of compact components of 
  $X_S$ by $\kappa(n-i,k)$ instead of $\kappa(n-i,k+1)$.
}\end{remark}

\section*{Acknowledgment}
We thank Alicia Dickenstein and Daniel Perucci whose comments inspired us to find the 
correct statement and proof of Theorem~\ref{Th:two}.



\newpage

\begin{thebibliography}{10}

\bibitem{blr}
{\sc R.~Benedetti, F.~Loeser, and J.-J. Risler}, {\em Bounding the number of
  connected components of a real algebraic set}, Discrete Comput. Geom., 6
  (1991), pp.~191--209.

\bibitem{Bihan}
{\sc F.~Bihan}, {\em {Polynomial systems supported on circuits and dessins
  d'enfants}}, 2005.
\newblock Journal of the London Mathematical Society, to appear.

\bibitem{BS}
{\sc F.~Bihan and F.~Sottile}, {\em {New fewnomial upper bounds from Gale dual
  polynomial systems}}, 2006.
\newblock Moscow Mathematical Journal, to appear. {\tt math.AG/0609544}.

\bibitem{GPS05}
{\sc G.-M. Greuel, G.~Pfister, and H.~Sch\"onemann}, {\em {\sc Singular} 3.0},
  {A Computer Algebra System for Polynomial Computations}, Centre for Computer
  Algebra, University of Kaiserslautern, 2005.
\newblock {\tt http://www.singular.uni-kl.de}.

\bibitem{Kh80}
{\sc A.~Khovanskii}, {\em A class of systems of transcendental equations},
  Dokl. Akad. Nauk. SSSR, 255 (1980), pp.~804--807.

\bibitem{Kh91}
\leavevmode\vrule height 2pt depth -1.6pt width 23pt, {\em Fewnomials}, Trans.
  of Math. Monographs, 88, AMS, 1991.

\bibitem{LRW03}
{\sc T.-Y. Li, J.~M. Rojas, and X.~Wang}, {\em Counting real connected
  components of trinomial curve intersections and {$m$}-nomial hypersurfaces},
  Discrete Comput. Geom., 30 (2003), pp.~379--414.

\bibitem{Na01}
{\sc D.~Napoletani}, {\em A power function approach to {K}ouchnirenko's
  conjecture}, in Symbolic computation: solving equations in algebra, geometry,
  and engineering (South Hadley, MA, 2000), vol.~286 of Contemp. Math., Amer.
  Math. Soc., Providence, RI, 2001, pp.~99--106.

\bibitem{Pe05}
{\sc D.~Perrucci}, {\em Some bounds for the number of components of real zero
  sets of sparse polynomials}, Discrete Comput. Geom., 34 (2005), pp.~475--495.

\bibitem{PW}
{\sc W.~Pohl and M.~Wenk}, {\em {\tt solve.lib}}, 2006.
\newblock {A {\sc Singular} 3.0 library for Complex Solving of Polynomial
  Systems}.

\end{thebibliography}
\end{document}